[1994/12/01]
\documentclass[draft]{article}
\title{Cuntz-like Algebras}
\author{Jean Renault}
\date{\it Dedicated to Professor Marc Rieffel on his 60th birthday}

\usepackage{amsmath,amsthm}

\chardef\bslash=`\\ 

\newtheorem{thm}{Theorem}[section]
\newtheorem{cor}[thm]{Corollary}
\newtheorem{lem}[thm]{Lemma}
\newtheorem{prop}[thm]{Proposition}

\theoremstyle{definition}
\newtheorem{defn}{Definition}[section]
\newtheorem{ex}{Example}[section]

\theoremstyle{remark}

\newcommand{\propref}[1]{Proposition~\ref{#1}}
\newcommand{\defnref}[1]{Definition~\ref{#1}}

\renewcommand{\labelenumi}{(\roman{enumi})}

\newcommand{\eval}[2][\right]{\relax
  \ifx#1\right\relax \left.\fi#2#1\rvert}

\begin{document}
\maketitle
\markboth{Cuntz-like Algebras}
{Cuntz-like Algebras}
\begin{abstract}
The usual crossed product construction which
associates to the homeomorphism $T$ of the locally compact space $X$
the C$^*$-algebra
$C^*(X,T)$ is extended to the case of a partial local homeomorphism $T$.
For example, the Cuntz-Krieger algebras are the C$^*$-algebras of the
one-sided Markov shifts. The generalizations of the Cuntz-Krieger algebras
(graph algebras, algebras $O_A$ where $A$ is an infinite matrix) which
have been introduced recently can also be described as C$^*$-algebras of
Markov chains with countably many states. This is useful to obtain such
properties of these algebras as  nuclearity, simplicity or
pure infiniteness. One also gives examples of strong Morita
equivalences arising from dynamical systems equivalences.\footnote{


Endomorphisms. Markov shifts. Morita equivalence.
\ }

\end{abstract} 
\renewcommand{\sectionmark}[1]{}

\section{Introduction.}

Let us recall (with no respect for history) two striking results
pertaining the rich interplay between ergodic theory and von Neumann
algebras (we refer the reader to the survey \cite{moo:relations} and the references
thereof for details; at that time, the theory was essentially complete
except the uniqueness of the hyperfinite $III_1$ factor). A
generalization \cite{fm:relations} of the Murray-von Neumann group measure construction
associates to each discrete measured equivalence relation (with
possibly a twist) a von Neumann algebra, which is injective (or
amenable) if and only if the measured equivalence relation is
amenable in the sense of Zimmer. This construction is particularly
well-behaved for amenable von Neumann algebras. First, each amenable
von Neumann algebra arises from a discrete measured equivalence relations.
Second, this discrete measured equivalence relation is unique up to
isomorphism. There is an important fact which underlies these developments and
also makes the connection with classical ergodic theory, namely: amenable equivalence
relations are singly generated \cite{cfw:amenable}.

Although the interplay between topological dynamics and C$^*$-algebras may be
more elusive, the study of transformation group C$^*$-algebras has proved for over
thirty years to be of great importance, both for the internal theory of
C$^*$-algebras and for its applications.  

In \cite{ar:examples}, we have proposed
\'etale essentially principal groupoids as a topological analogue of discrete
measured equivalence relations; just as discrete equivalence relations arise from
group actions, these groupoids arise as groupoids of germs of pseudogroups. Via the
groupoid construction of \cite{ren:approach} (or the equivalent localization
construction of
\cite{kum:localizations}) , one obtains a large class of C$^*$-algebras, which are
nuclear if and only if the groupoid is amenable in the sense of \cite{dr:amenable}. We
shall study here a subclass of these algebras, namely those  arising from singly
generated pseudogroups (a precise definition will be given below). We view the groupoids
of germs of singly generated pseudogroups as a (sort of) topological analogue of singly
generated equivalence relations;  in particular, shifts and shift spaces will provide
our basic examples.

In the first section, we shall introduce this class of ``Cuntz-like'' algebras and show
that they share some of the features of the classical Cuntz algebras $O_n$. In the
second section, we consider some groupoid equivalences, in particular those arising from
shift equivalences. The results presented there are well known. In the third section, we
shall describe as groupoid C$^*$-algebras the examples of the Cuntz-Krieger algebras for
infinite matrices constructed by R.~Exel and M.~Laca in \cite{el:infinitematrices}.
These examples were a strong motivation for our definitions. 
 
This work, which is an extension of a talk given
at the 17th Conference on Operator Theory at Timisoara in June 98, illustrates the use
of  some groupoid techniques in the study of Cuntz-Krieger algebras. It only covers a
limited part of the rich domain of the Cuntz-Krieger algebras and their generalizations.

I heartily thank Claire Anantharaman-Delaroche, Marcelo Laca and Georges Skandalis for
valuable comments. 

\section{Singly generated pseudogroups and their C$^*$-algebras.}

We first recall some definitions from \cite{ar:examples}.

\begin{defn} 

Let $X$ be a topological space. A partial homeomorphism $S$
of $X$ is a homeomorphism from an open set $dom(S)$ of $X$ onto an open set
$ran(S)$ of $X$. A pseudogroup $\cal G$ on $X$ is a set of partial homeomorphisms
of $X$ such that
\begin{enumerate}
\item the identity map $1_X$ belongs to $\cal G$,
\item $S,T\in {\cal G}\Rightarrow ST=S\circ T\in {\cal G}$,
\item  $S\in {\cal G}\Rightarrow S^{-1}\in {\cal G}$.

The pseudogroup
$\cal G$ is said to be full if 
\item for every open set $U$ of $X$, $S\in {\cal G}\Rightarrow S_{|U}=S\ 1_U\in
{\cal G}$,
\item  each partial homeomorphism $S$ which belongs locally to $\cal G$, i.e.
every $x\in dom(S)$ has an open neighborhood $U_x$ such that the restriction of $S$
to $U_x$ belongs to $\cal G$, does belong to $\cal G$.
\end{enumerate}

\end{defn}

The conditions $(iv)$ and $(v)$ are usually included in  the definition of a
pseudogroup (cf. \cite{hae:pseudogroups}). By analogy with the notion of full group in
ergodic theory, we prefer to state them apart. We remark that every pseudogroup $\cal G$
is contained in a smallest full pseudogroup $\overline{\cal G}$. We say that
$\overline{\cal G}$ is the full pseudogroup of $\cal G$. There are several groupoids,
related to the semidirect product of $X$ by $\cal G$, which can be attached to the
pseudogroup $\cal G$. We are chiefly interested in its groupoid
of germs, which depend only on the full pseudogroup. Let us recall its definition.

\begin{defn}
We say that two partial
homeomorphisms $S,T$ have the same germ at
$x\in dom(S)\cap dom(T)$ if $S$ and $T$ agree on an open neighborhood of $x$. The
equivalence class is called the germ at $x$ and denoted by $[Sx,S,x]$.  The groupoid
of germs
$G=Germ(X,{\cal G})$ of the pseudogroup
$\cal G$ on the topological space
$X$ is the set of germs of $\cal G$:
 $$G=\{[x,S,y], S\in{\cal G}, y\in dom(S), x=Sy\}$$
endowed with the groupoid structure defined by the maps $r,s:G\rightarrow X$ such
that
$r[x,S,y]=x, s[x,S,y]=y$, the product $[x,S,y][y,T,z]=[x,ST,z]$ and the inverse
$[x,S,y]^{-1}= [y,S^{-1},x]$ and the topology of germs with 
$${\cal U}(S;U)=\{[x,S,y]: y\in U\}$$
where $S\in{\cal G}$ and $U$ is an open subset of
$dom(S)$ as basic open sets.
\end{defn}
If $X$ is a Baire space, then the groupoid of germs $Germ(X,{\cal G})$ is essentially
principal in the sense that the set of points with trivial isotropy is dense (cf.
Proposition~2.1 of
\cite{ar:examples}).

\begin{defn} We define a {\it singly generated dynamical system} (SGDS) as a
pair $(X,T)$ where $X$ is a topological space and
$T$ is a local homeomorphism from an open subset $dom(T)$ of $X$ onto an open
subset $ran(T)$ of $X$. We denote by ${\cal G}(X,T)$ the full pseudogroup generated by
the restrictions $T_{|U}$, where $U$ is an open subset of
$X$ on which $T$ is injective.
\end{defn}
\begin{lem} 
Let $(X,T)$ be a SGDS.
\begin{enumerate}
\item A partial homeomorphism $S$ belongs to ${\cal
G} (X,T)$ iff it is locally of the form $(T^m_{|U})^{-1}T^n_{|V}$, where $m,n\in{\bf N}$,
$U$ is an open set on which $T^m$ is injective and $V$ is an open set on which
$T^n$ is injective.
\item Let $x\in X$. Suppose that  $(T^m_{|U})^{-1}T^n_{|V}$ and
$(T^p_{|W})^{-1}T^q_{|Y}$ are two partial homeomorphisms as in $(i)$ having $x$ in their
domains and sending it into the same element. If $m-n=p-q$, then $(T^m_{|U})^{-1}T^n_{|V}$ and
$(T^p_{|W})^{-1}T^q_{|Y}$ have the same germ at $x$.
\end{enumerate}
\end{lem}

\begin{proof} 
$(i)$ It is clear that a partial homeomorphism of the above form
$(T^m_{|U})^{-1}T^n_{|V}$ is in ${\cal G}(X,T)$ and, since ${\cal G}(X,T)$ is full, this is
still true for a partial homeomorphism locally of that form. On the other hand, the
inverse of $(T^m_{|U})^{-1}T^n_{|V}$, which is $(T^n_{|V})^{-1}T^m_{|U}$, is of the
same form. The product  of $(T^m_{|U})^{-1}T^n_{|V}$ and $(T^p_{|W})^{-1}T^q_{|Y}$ is
also of the same form: if $n\ge p$, it is of the form 
$(T^m_{|U})^{-1}T^{n-p+q}_{|Z}$, where $Z$ is an open set on which $T^{n-p+q}$ is
injective and if $n<p$, it is of the form  $(T^{p-n+m}_{|Z})^{-1}T^q_{|Y}$, where $Z$ is
an open set on which $T^{p-n+m}$ is injective. This implies that
the set of partial homeomorphisms locally of the above form is a full pseudogroup,
and therefore is ${\cal G}(X,T)$.

$(ii)$ We may assume that $m-p=n-q=k\ge 0$, $W=U, Y=V$ and $T^pU=T^qV$
Then 
$$(T^m_{|U})^{-1}T^n_{|V}=(T^p_{|U})^{-1}(T^k_{|T^pU})^{-1}T^k_{|T^qV}T^q_{|V}
=(T^p_{|U})^{-1}T^q_{|V}.$$
\end{proof}

 We denote by $Germ(X,T)$ the groupoid of germs of ${\cal G}(X,T)$. We also consider
another groupoid attached to $(X,T)$, the semidirect product groupoid:

\begin{defn} 
Let $(X,T)$ be a SGDS. Its semidirect product groupoid is
$$G(X,T)=\{(x,m-n,y): m,n\in {\bf N},  x\in dom(T^m), y\in dom(T^n), T^mx=T^ny\}$$
with the groupoid structure induced by the product structure
of the trivial groupoid $X\times X$ and of the group ${\bf Z}$ and the topology defined
by the basic open sets
$${\cal U}(U;m,n;V)=\{(x,m-n,y): (x,y)\in U\times V,\ T^m(x)=T^n(y)\}$$ 
where $U$  [resp.
$V$] is an open subset of the domain of $T^m$ [resp. $T^n$] on
which $T^m$ [resp. $T^n$] is injective.
\end{defn}
 According to the above lemma, there is a
map $\pi$ from $G(X,T)$ onto $Germ(X,T)$ which sends $(x,m-n,y)$ into the germ 
$[x,(T^m_{|U})^{-1}T^n_{|V},y]$, where $U$ is an open neighborhood of $x$ on which $T^m$
is injective and
$V$ is an open neighborhood of $y$ on which
$T^n$ is injective. This map is continuous and is a groupoid homomorphism. Let us see
when it is an isomorphism.

\begin{defn} We shall say that a SGDS $(X,T)$ is essentially free if for
every pair of distinct integers $(m,n)$, there is no nonempty open set on which
$T^m$ and $T^n$ agree.
\end{defn}

\begin{lem} Let $(X,T)$ be an essentially free SGDS. Then,
\begin{enumerate}
\item If $(T^m_{|U})^{-1}T^n_{|V}$ and $(T^p_{|W})^{-1}T^n_{|Y}$, where
$m,n,p,q\in{\bf N}$ and $U,V,W,Y$ are open sets such that 
$T^m_{|U},T^n_{|V},T^p_{|W},T^n_{|Y}$ are injective, have the
same germ at
$x$, then
$m-n=p-q$.
\item The map $c:Germ(X,T)\rightarrow {\bf Z}$  such that
$c[(T^m_{|U})^{-1}T^n_{|V}x, (T^m_{|U})^{-1}T^n_{|V},x]=m-n$ is a continuous
homomorphism.
\end{enumerate}
\end{lem}

\begin{proof} 
$(i)$ By assumption, the relation $T^{n+p}y=T^{m+q}y$ holds on a
nonempty open set. The essential freeness implies that $n+p=m+q$.

$(ii)$ We have seen in the previous proposition that the product  of
$(T^m_{|U})^{-1}T^n_{|V}$ and $(T^p_{|W})^{-1}T^q_{|Y}$ is  of the  form 
$(T^m_{|U})^{-1}T^{n-p+q}_{|Z}$ or  $(T^{p-n+m}_{|Z})^{-1}T^q_{|Y}$ and that the inverse
of  $(T^m_{|U})^{-1}T^n_{|V}$ is $(T^n_{|V})^{-1}T^m_{|U}$. This shows that $c$ is a
homomorphism. By construction, $c$ is locally constant.
\end{proof}

\begin{prop} Let $(X,T)$ be a SGDS. Then, 
$(X,T)$ is essentially free if and only if the above surjection 
$\pi:G(X,T)\rightarrow Germ(X,T)$ is an isomorphism.
\end{prop}

\begin{proof}
If $(X,T)$ is essentially free, we may define the map $\rho$ from $Germ(X,T)$ to
$G(X,T)$ by $\rho[y,S,x]=(y,c(y,S,x),x)$. This
map is an inverse for $\pi$. 
 
Conversely, suppose that $\pi$ is an isomorphism. Let $m,n\in{\bf N}$ and $U$ be a
nonempty open set such that $T_{|U}^m=T_{|U}^n$. Then for $x\in U$, the germ
$\pi(x,m-n,x)=[x,(T^m_{|U})^{-1}T^n_{|U},x]$ is a unit. This implies that $m=n$. 
\end{proof}

From now on, we shall assume that the space $X$ of our  SGDS $(X,T)$  is Hausdorff,
second countable and locally compact. Then
$G(X,T)$ is a Hausdorff locally compact \'etale groupoid and we can construct
 its  C$^*$-algebra
$C^*(X,T)=C^*(G(X,T))$ (we shall soon see that we do not have to
distinguish the full and the reduced C$^*$-algebras). 

Here is an elementary example of a SGDS. We let $X$ be the one-point compactification
of
${\bf N}\cup\{\infty\}$ and $T:{\bf N}\rightarrow {\bf N}$ be the translation $i\mapsto
i+1$.  Then $(X,T)$ is essentially principal and $G(X,T)=Germ(X,T)=({\bf N}\times{\bf
N})\cup\{(\infty,\infty)\}$ where the sets 
$\{(i,i):i\ge i_0\}$ form a fundamental system of neighborhoods of $\{(\infty,\infty)\}$.
Its C$^*$-algebra is the algebra obtained by ajoining a unit to the algebra of compact
operators on $l^2({\bf N})$.
 
An other example of SGDS is provided by the full one-sided shift, where 
$X=\{1,\ldots,d\}^{\bf N}$,  $d$ is an integer, and $Tx_i=x_{i+1}$. Then $C^*(X,T)$ is
the Cuntz algebra $O_d$ (cf. Section~III.2 of \cite{ren:approach}). A generalization
of this example will be discussed in the third section. Some of the features of the Cuntz
algebras are kept in the general situation. The crossed product nature of
$C^*(X,T)$ is revealed by the fundamental homomorphism $c:G(X,T)\rightarrow{\bf Z}$,
which induces the dual action (cf. \cite{ren:approach}, Section~II.5) 
$\alpha$ of ${\bf T}=\hat{\bf Z}$ on $C^*(X,T)$ according to 
$\alpha_t(f)(x,k,y)=e^{ikt}f(x,k,y)$ for $f\in C_c(G(X,T))$. The kernel of $c$ will be
denoted
$R(X,T)=c^{-1}(0)$. It is an open and closed subgroupoid of $G(X,T)$ and it has no isotropy.
It is reduced to the unit space $X$ when $T$ is a global homeomorphism of $X$ onto $X$.
Its C$^*$-algebra $C^*(R(X,T))$ is the fixed point algebra $C^*(X,T)^\alpha$. In the case
of the full one-sided shift as above, $C^*(R(X,T))$ is the UHF algebra $UHF(d^\infty)$.

\begin{prop} Let $(X,T)$ be a SGDS. Then,
\begin{enumerate}
\item $G(X,T)$ is amenable;
\item the full and reduced C$^*$-algebras coincide;
\item the C$^*$-algebra $C^*(X,T)$ is nuclear.
\end{enumerate}
\end{prop}

\begin{proof}
$(i)$ We will check measurewise amenability ( according to \cite{dr:amenable}, 3.3.7, it
is equivalent to topological amenability for \'etale groupoids). We first show that
$R(X,T)=c^{-1}(0)$ is amenable. It is the  increasing union of $R_N=\{(x,y)\in X\times
X:\exists n\le N: x,y\in dom(T^n), T^nx=T^ny\}$. According to Section~5.f of
\cite{dr:amenable}, it suffices to show that $R_N$ is amenable. This is a Borel
equivalence relation with countable equivalence classes. There is a countable family of
open  sets $U_i$ which covers 
$X$ and such that, for each $n\in\{0,\ldots, N\}$, the restrictions $T^n_{|U_i}$ are
one-to-one. Therefore, the equivalence relation $R_N$ is countably separated, its
quotient space
$X/R_N$ is analytic . This implies that $R_N$ is a proper Borel
groupoid (\cite{dr:amenable}, 2.1.2), hence it is amenable. We can now apply a result on
the amenability of an extension (\cite{dr:amenable}; 5.2.13) to conclude that $G(X,T)$
is amenable. Indeed, we may write $X$ as the disjoint union of the invariant Borel
subsets $Y$ and
$Z$, where $Z$ is the intersection of the domains of the $T^n$'s and $Y$ is its
complement. The homomorphism
$c:G(X,T)\rightarrow {\bf Z}$ is strongly surjective in the sense given there on $Z$. On
the other hand, the reduction of $G(X,T)$ to $Y$ is a proper principal groupoid having as
quotient space the complement $U^c$ of the domain $U$.

$(ii)$ This is a well-known property of the C$^*$-algebra of an amenable groupoid (see
e.g. \cite{dr:amenable}; 6.1.5).

$(iii)$  This is also a well-known property of the C$^*$-algebra of an amenable groupoid
(see e.g. \cite{dr:amenable}; 6.2.14).
\end{proof}

According to \cite{ren:ideal}, since $G(X,T)$ is Hausdorff and amenable, the ideal
structure of
$C^*(X,T)$ is very much related to the structure of the open invariant subsets of $X$.
Note that invariance with respect to a groupoid of germs
$G$ of
a pseudogroup $\cal G$ simply means invariance under $\cal G$. A sufficient condition
(\cite{ren:ideal}, 4.9) for an isomorphism between these two structures is that for every
closed invariant subset of $X$, the reduced groupoid $G(X,T)_{|F}$ has a dense set of
points without isotropy (a word of caution has to be given here: the reduced
groupoid $G_{|F}$ is usually distinct from the groupoid of germs of the reduced
pseudogroup ${\cal G}_{|F}$). In particular, we have the following criterion for the
simplicity of $C^*(X,T)$:

\begin{prop}\label{simplicity} Let $(X,T)$ be an essentially free SGDS. Assume that for
every nonempty open set $U\subset X$ and every $x\in X$, there exist $m,n\in{\bf N}$
such that $x\in dom(T^n)$ and $T^nx\in T^mU$. Then $C^*(X,T)$ is simple.
\end{prop}
We can also express in our particular case the locally contracting property of 
Section 2.~4 of \cite{del:purelyinfinite} which ensures that the C$^*$-algebra is purely
infinite.

\begin{prop}\label{pure infiniteness} Let $(X,T)$ be an essentially free SGDS. Assume
that for every nonempty open set $U\subset X$, there exist an open set $V\subset U$ and
$m,n\in{\bf N}$ such that $T^n(V)$ is strictly contained in $T^m(V)$. Then $C^*(X,T)$
is purely infinite.
\end{prop}

\section{Groupoid equivalence and shift equivalence.}

\subsection{Equivalence of pseudogroups}

\begin{defn}
Let $(X,\cal G)$ and $(Y,\cal H)$ be two full pseudogroups. An
isomorphism from
$\cal G$ to $\cal H$ is  a homeomorphism  $\varphi:X\rightarrow Y$ such that
$\varphi {\cal G}\varphi^{-1}={\cal H}$.
\end{defn}

It is useful to introduce a notion weaker than isomorphism of pseudogroups.

\begin{defn} Let $(X,\cal G)$ be a full pseudogroup and let $A,B$ be open
subsets of $X$. We denote by
${\cal G}_A^B$ the set of elements $S\in {\cal G}$ which have their domains  contained
in $B$ and their ranges contained in $A$. In particular, ${\cal G}_{|A}={\cal
G}_A^A$ is called the reduction of $\cal G$ to $A$. We say that $A\subset X$ is full
(with respect ot $\cal G$) if it meets every orbit under $\cal G$.
\end{defn}

\begin{defn}(Cf. \cite{hae:pseudogroups}) Let $(X,\cal G)$ and $(Y,\cal H)$ be two full
pseudogroups. An equivalence from $\cal G$ to $\cal H$ is a maximal collection $\cal Z$
of partial homeomorphisms from $X$ to $Y$ such that
\begin{enumerate}
\item The domains of the elements of $\cal Z$ form an open cover of $X$ and the
ranges of the elements of $\cal Z$ form an open cover of $Y$.
\item $R\in {\cal G},\  \varphi\in {\cal Z}, \  T\in {\cal H}\Longrightarrow
T\varphi R\in {\cal Z}$.
\item $R\in {\cal G},\  \varphi,\psi\in {\cal Z}, \  T\in {\cal
H}\Longrightarrow
\psi R\varphi^{-1}\in {\cal G},\  \psi^{-1} T\varphi\in {\cal H} $.
\end{enumerate}
\end{defn}

 Let $(X,\cal G)$ and $(Y,\cal H)$ be two full pseudogroups and let $\cal Z$ be an
equivalence from $\cal G$ to $\cal H$. Let $\cal I$ be the full pseudogroup on the
disjoint union $Z=X\sqcup Y$ generated by $\cal G,H,Z$. Then $X$ and $Y$ are full with
respect to $\cal I$ and ${\cal G}={\cal I}_{|X}, {\cal H}={\cal I}_{|Y}$. Conversely, let $(Z,
{\cal I})$ be a full pseudogroup and $X,Y$ two full open subsets of $Z$. Then ${\cal
I}_Y^X$ is an equivalence from ${\cal I}_{|X}$ onto ${\cal I}_{|Y}$.

Let $(X,\cal G)$ and $(Y,\cal H)$ be two pseudogroups. Every collection
${\cal Z}$ of partial homeomorphisms from $Y$ to $X$ such that
\begin{enumerate}
\item The union of the domains of the elements of ${\cal Z}$ meets every orbit 
under
$\cal G$ and the union of the ranges of the elements of $\cal Z$ meets every orbit 
under
$\cal H$.
\item  $S\in {\cal H},\ \varphi,\psi \in {\cal Z},\  T\in {\cal
G}\Longrightarrow
\varphi T\psi^{-1}\in \overline{\cal H},\  \varphi^{-1} S\psi\in \overline{\cal G}$
\end{enumerate}
can be completed in a unique way into an equivalence
$\overline{\cal Z}$ between
$\overline{\cal H}$ and $\overline{\cal G}$.
Explicitly,
 a partial homeomorphism $\varphi$ from $Y$ to $X$ belongs to $\overline{\cal Z}$ iff
for every point $y$ of its domain, there exist an open neighborhood $U$ of $y$,
$S\in {\cal H}, \psi \in {\cal Z}, T\in {\cal G}$  such that
$\varphi_{|U}=S\psi T_{|U}$. We then say that the collection  ${\cal Z}$ generates
 $\overline{\cal Z}$.

The proof of the following proposition is straightforward.

\begin{prop}
Let $(X,\cal G)$ and $(Y,\cal H)$ be two full pseudogroups and
let
$G$ and $H$ be their groupoids of germs.
\begin{enumerate} 
\item An isomorphism $\varphi:X\rightarrow Y$ from $\cal G$ onto $\cal H$
implements the groupoid isomorphism $\varphi:G\rightarrow H$ such that
$\varphi[x,S,y]= [\varphi(x),\varphi S\varphi^{-1},\varphi(y)]$.
\item An equivalence $\cal Z$ from $\cal G$ onto $\cal H$
implements the groupoid equivalence 
$Z$ from $G$ onto $H$, where $Z$ is the space of germs of $\cal Z$.
\end{enumerate}
\end{prop}

\begin{ex}{\em The irrational rotations algebras.} (cf. \cite{rie:rotations} and
\cite{kum:localizations})  Let $a$ be an irrational number. Then the  pseudogroup
$\cal G$ on $X={\bf R}$ generated by the translations $x\mapsto x+1$ and $x\mapsto x+a$
is equivalent to the pseudogroup $\cal H$ on $Y={\bf R}/{\bf Z}$ generated by the
translation  $y\mapsto y+a$. Indeed the restriction of the quotient map
$X\rightarrow Y$ to a nonempty open interval of length strictly less than one generates
an equivalence from $\cal G$ onto $\cal H$. Therefore, if $a$ and $b$ are irrational
numbers such that ${\bf Z}+a{\bf Z}= {\bf Z}+b{\bf Z}$, then the groupoids of the
corresponding irrational rotations are equivalent and their C$^*$-algebras are strongly
Morita equivalent.
\end{ex}

\subsection{Strong shift equivalence.}

We shall limit ourselves in this section to SGDS $(X,S)$such that\hfill\break
$dom(S)=ran(S)=X$. We will investigate further the structure of the semidirect product
$G(X,S)$. We have introduced the fundamental homomorphism (or cocycle)
$c:G(X,S)\rightarrow {\bf Z}$ such that $c(x,m-n,y)=m-n$ and its kernel $R(X,S)$.
Another important piece of structure is the fundamental endomorphism of $G(X,S)$,
which is still denoted by $S$ and which is defined by $S(x,m-n,y)=(S(x),m-n,S(y))$. This
is indeed a groupoid homomorphism of $G(X,S)$ into itself; it is surjective since
$S:X\rightarrow X$ is assumed to be surjective; it preserves the fundamental cocycle in
the sense that $c\circ S=c$ and induces an endomorphism of $R(X,S)$ onto itself. As an
endomorphism of $G(X,S)$, it is equivalent to the identity since we have
$S(\gamma)=(b\circ r(\gamma))^{-1}\gamma b\circ s(\gamma)$ where $b(x)=
(x,1,Sx)$. Its graph is a groupoid equivalence isomorphic to the identity equivalence.
However, in general, it is no longer trivial as an endomorphism of $R(X,S)$.

\begin{lem} 
The graph $Z(S)$ of the fundamental homomorphism is a groupoid
equivalence from
$R(X,S)$ onto itself.
\end{lem}

\begin{proof} 
This graph is the space 
\[Z(S)=\{(x,y)\in X\times X:\exists m: S^{m+1}(x)=S^m(y)\}\]
endowed with the natural left and right actions of $R(X,S)$. One can check directly the
axioms of \cite{mrw:equivalence}.
\end{proof}

Given the SGDS $(X,S)$, we form the SGDS $(\tilde X,\tilde S)$
where $\tilde X$ is the projective limit of the projective system 
\[\ldots X_2\xrightarrow{\pi_{1,2}} X_1\xrightarrow{\pi_{0,1}} X_0\]
where $X_i=X$ and $\pi_{i,i+1}=S$ for all $i=0,1,\ldots$ and $\tilde S$ is the
homeomorphism induced by the projective system morphism $S_i=S:X_i\rightarrow
X_i$. Then $\tilde S$ is invertible and its inverse is defined
by the projective system morphism $T_i=id:X_{i+1}\rightarrow
X_i$. We denote by $\pi=\pi_0$ the projection $\tilde X\rightarrow X=X_0$ and call
$(\tilde X,\tilde S)$ the invertible extension of $(X,S)$.

\begin{defn}
Let $(X,S)$ and $(Y,T)$ be two SGDS. A (strong) shift equivalence
of lag
$k\in{\bf N}$, consists of a pair of continuous maps
$\varphi:X\rightarrow Y$ and $\psi:Y\rightarrow X$ such that
\[\begin{array}{cc}
T\varphi=\varphi S\ ;&\psi T=S\psi\\
\psi\varphi=S^k\ ;&\varphi\psi=T^k.
\end{array}\]
\end{defn}

The definition implies that $\varphi$ and $\psi$ are surjective and
local homeomorphisms. When $S$ and $T$ are invertible,  $\varphi$ and $\psi$ are
simply conjugacies. When $S$ and $T$ are not invertible, 
the strong shift equivalence $(\varphi,\psi)$ from $(X,S)$ to $(Y,T)$ induces a pair of
morphisms of the associated projective systems, hence a pair of continuous maps
$(\tilde\varphi,\tilde\psi)$ from $(\tilde X,\tilde S)$ to $(\tilde Y,\tilde T)$ which
is also a strong shift equivalence with the same lag.
Therefore,
$\tilde\varphi$ and $\tilde\psi$ are conjugacies.
 
Let us look at the strong shift equivalence from the groupoid point of view. The
condition
$T\varphi=\varphi S$ implies that
$\varphi$ induces a groupoid homomorphism, still denoted by $\varphi$, from $G(X,S)$
onto
$G(Y,T)$ such that $\varphi(x,m-n,y)=(\varphi(x),m-n,\varphi(y))$.  Since it preserves
the fundamental cocycles, it also defines  a groupoid homomorphism from $R(X,S)$
onto $R(Y,T)$. Its graph is a groupoid equivalence which intertwines the fundamental
equivalences $Z(S)$ and $Z(T)$. 
 
At the level of the C$^*$-algebras, a strong shift equivalence implements a
strong Morita equivalence between $C^*(X,S)$ and $C^*(Y,T)$ which is trivial in K-theory
and a strong Morita
equivalence between $C^*(R(X,S))$ and $C^*(R(Y,T))$ which induces an isomorphism of
their K-groups which intertwines the automorphisms induced by $S$ and $T$.

\begin{ex}
{\em Topological Markov Shifts.}
We first recall from \cite{lm:intro} the definition of the edge shift associated with a
graph. Let
$\Gamma$ be a directed graph with finite vertex set
$I$ and finite edge set
$E$. The  adjacency matrix of the graph is the nonnegative integer valued matrix $A$
whose element $A(i,j)$ is the number of edges from vertex $i$ to vertex $j$. Up to
isomorphism $\Gamma$ depends only on $A$ and we write $\Gamma=\Gamma_A$.  We
denote by
$X_A$ the space of one-sided infinite paths
$e_0e_1e_2\ldots$, where
$e_n\in E$ and $r(e_n)=s(e_{n+1})$ endowed with the product topology and by $T_A$ the
one-sided shift $(T_Ae)_n=e_{n+1}$ on $X_A$.

Recall that an elementary equivalence  between two 
nonnegative integral square matrices
$A$ and
$B$  is a pair $(R,S)$ of nonnegative integral matrices such that $A=RS$ and $B=SR$.
A strong shift equivalence between
$A$ and
$B$ is a finite sequence of 
nonnegative integral square matrices
$A_0=A,A_1,\ldots,A_k=B$ and elementary equivalences between $A_{i-1}$ and $A_i$,
for $i=1,\ldots,k$.  A strong shift equivalence between the matrices $A$ and $B$ gives
a strong shift equivalence between the SGDS $(X_A,T_A)$ and $(X_B,T_B)$. It suffices
to consider the case of an elementary equivalence $(R,S)$. One draws $R(i,j)$ edges
from the vertex $i$ of $\Gamma_A$ to the vertex $j$ of $\Gamma_B$ and call $E_R$
this set of edges. One constructs similarly a set of edges $E_S$. The condition $A=RS$
gives the existence of a bijection $\alpha$ from $E_A$ onto the set $E_R*E_S$ of paths
$rs$ with $r\in E_R,s\in E_S$ preserving the initial and terminal vertices. Similarly,
one chooses a bijection $\beta$ from $E_B$ onto $E_S*E_R$. One defines
$\varphi:X_A\rightarrow X_B$ by $\varphi(a_0a_1\ldots)=b_0b_1\ldots$ where
$\alpha(a_0)=r_0s_0,\alpha(a_1)=r_1s_1,\ldots$ and
$\beta(b_0)=s_0r_1,\beta(b_1)=s_1r_2,\ldots$. One defines similarly
$\psi:X_B\rightarrow X_A$. One checks that $\varphi$ and $\psi$ are continuous and
that $T_A=\psi\varphi$ and $T_B=\varphi\psi$. Therefore, according to the above, a
strong shift equivalence between the matrices $A$ and $B$ gives a groupoid equivalence
between the principal groupoids $R(X_A,T_A)$ and $R(X_B,T_B)$ which intertwines the
fundamental equivalences $Z(T_A)$ and $Z(T_B)$. In this example, 
$C^*(R(X_A,T_A))$ is an AF-algebra having for dimension group the inductive limit
\[G(A)=\lim_{\rightarrow}({\bf Z}^{n_A}\xrightarrow{A} {\bf
Z}^{n_A}\xrightarrow{A}\ldots)\] and
$T_A$ induces the shift automorphism $\tau_A$.
Thus we have a realization at the groupoid level of the shift preserving isomorphisms of
the dimension groups $G(A)$ and $G(B)$ induced by the strong shift equivalence.

On the other hand, recall that a shift equivalence of lag $k$  between two 
nonnegative integral square matrices
$A$ and
$B$  is a pair $(R,S)$ of nonnegative integral matrices such that
\[\begin{array}{cc}
AR=RB\ ;&SA=BS\\
RS=A^k\ ;&SR=B^k.
\end{array}\]
In particular $(R,S)$ is an elementary equivalence between $A^k$ and $B^k$. Note
that\hfill\break
$(X_{A^k},T_{A^k})$ is conjugate to $(X_A,T_A^k)$ and that $R(X_A,T_A^k)=R(X_A,T_A)$.
Thus, if $A$ and $B$ are only shift equivalent, we can still construct a groupoid
equivalence between $R(X_A,T_A)$ and $R(X_B,T_B)$. Because of the other two
relations, it induces a shift preserving isomorphisms of
the dimension groups $G(A)$ and $G(B)$ (in fact $A$ and $B$ are shift equivalent if and
only if $(G(A),G(A)^+,\tau_A)$ and $(G(B),G(B)^+,\tau_B)$ are isomorphic
\cite{kri:dimensions}. This equivalence intertwines the  equivalences
$Z(T_A^k)$ and $Z(T_B^k)$ but not necessarily the equivalences
$Z(T_A)$ and $Z(T_B)$.
\end{ex}

\section{Graphs and Cuntz-Krieger algebras.}

We describe in this section a class of SGDS $(X,T)$ which appears implicitly in the work
\cite{el:infinitematrices}
of R.~Exel and M.~Laca. The Cuntz-Krieger
algebras for infinite matrices which they construct are the C$^*$-algebras of these
SGDS.

\begin{defn}\label{partition} 
Let $X$ be a compact totally
disconnected space, $U,V$ two open subsets and let $T:U\rightarrow V$ be a local
homeomorphism from $U$ onto $V$.
A Markov partition for $(X,T)$ is  a partition of
$U$ by a family
$\{U_i, i\in I\}$  of nonempty pairwise disjoint compact open subsets such that
\begin{enumerate}
\item the restriction $T_i=T_{|U_i}$ is a homeomorphism from
$U_i$ onto a compact open subset $V_i =T(U_i)$;
\item for all $(i,j)\in I\times I$, either $U_i\subset V_j$ or $U_i\cap V_j=
\emptyset$;
\item the Boolean algebra ${\cal B}_0$ generated by $\{X, U_i, V_i, i\in I\}$ is
a generator in the sense that $\bigvee_{n=0}^\infty T^{-n}{\cal B}_0$ is the family of
all compact open subsets of $X$.
\end{enumerate}
\end{defn}

\begin{defn}\label{markov} A Markov shift is a SGDS $(X,T)$ which admits a Markov
partition and which has a dense domain.
\end{defn}

The analysis of a SGDS $(X,T)$ admitting a Markov partition $\{U_i\}, i\in I$ can be done
entirely from the subset $A\subset I\times I$ defined by 
$(j,i)\in A \Leftrightarrow V_j\supset U_i$ and the subset 
${\cal J}\subset {\cal P}(I)$ defined by
$$J\in {\cal J}\Leftrightarrow
\cap_JV_j\cap_{I\setminus J}V^c_k\cap U^c\not=\emptyset .$$
Because of the condition $(iii)$ of Definition~\ref{partition}, this
intersection contains at most one point, which will be denoted by
$x_J$ for $J\in {\cal J}$.

Let us first determine the spectrum $X_0$ of the Boolean algebra ${\cal B}_0$ generated
by
$\{X, U_i, V_i, i\in I\}$.
Given $x\in X$, either $x\in
U$ and then there exists a unique $i$, which we denote by
$i=i_0(x)$, such that
$x\in U_i$   or
$x\in U^c$ and then we set $I(x)=\{j\in I: x\in V_j\}$. Note that $I(x)$
belongs to
$\cal J$ and that $x\mapsto I(x)$ is a bijection from $U^c$ onto $\cal J$, whose inverse
map is $J\mapsto x_J$. Putting these two maps together, we get a bijection, denoted by
$\sigma_0$, from
$X_0$ onto the disjoint union
$I\sqcup{\cal J}$. It remains to describe its topology. On $I$, it is the discrete topology.
On the other hand,
$\cal J$ is a closed subset of
${\cal P}(I)={\bf 2}^I$ endowed with the product
topology and the map: $\sigma_0$ is a
homeomorphism from $U^c$ onto
${\cal J}$. We introduce the map
$J:I\rightarrow {\cal P}(I)$ defined by $J(i)=\{j\in I: V_j\supset U_i\}$ and note that
a sequence $(x_\lambda)$ in $U$ converges to $x\in U^c$ iff $(J(i_0(x_\lambda)))$
converges to $\sigma_0(x)$. This shows that $\cal J$ contains the cluster points of the
map $J$ and that the topology of
$I\sqcup{\cal J}$ turning the symbol map $\sigma_0$ into a homeomorphism is the
topology induced by the compact space $I^+\times {\cal P}(I)$, where $I^+$ is
the one-point compactification of $I$ and $I^+\times {\cal P}(I)$ has the
product topology and where $I\sqcup_J{\cal J}$ is
identified with the closed subset $Graph(J)\cup\{\infty\}\times{\cal J}$. We denote by
$\tilde X_0=I\sqcup_J{\cal J}$ this topological space.
 
A similar analysis provides the spectrum $X_n$ of the
Boolean algebra ${\cal B}_n=\bigvee_{i=0}^n T^{-i}{\cal B}_0$: every $x\in X$ defines a
character of ${\cal B}_n$: one considers its orbit $(x,Tx,
T^2x,\ldots)$. If $x\in T^{-n}(U)$, then $x$ is coded by the sequence 
$\sigma_n(x)=(i_0(x),\ldots,i_n(x))$
where $i_k(x)=i_0(T^kx)$. If not, there exists an exit time $\tau=\tau(x)\le n$ such that
$T^\tau x\notin U$; then $x$ is coded by the sequence 
$\sigma_n(x)=(i_0(x),\ldots,i_{\tau-1}(x); I(T^\tau x))$. Thus we obtain a bijection
$\sigma_n$ from
$X_n$ onto the disjoint union $A^{(n)}\sqcup Y_n\ldots\sqcup Y_1\sqcup{\cal J}$,
where $A^{(m)}=\{(i_0,\ldots,i_m)\in I^{m+1}: (i_k,i_{k+1})\in A$ and
$Y_r=\{(i_0,\ldots,i_{r-1};J)\in A^{(r-1)}\times{\cal J}: i_{r-1}\in J\}$. Let us describe
the topology of this disjoint union. As before, $A^{(n)}$ has the discrete topology and
$Y_r$ has the topology of $A^{(r-1)}\times{\cal J}$. We have a sequence of continuous
maps
$$A^{(n)}\xrightarrow {J_n}A^{(n-1)}\times{\cal P}(I)\rightarrow\ldots\rightarrow
A^{(0)}\times{\cal P}(I)\xrightarrow {J_0}{\cal P}(I)$$
where $J_m(i_0,\ldots,i_m;J)=(i_0,\ldots,i_{m-1};J(i_m))$. We view this disjoint union
as a closed subset of ${A^{(n)}}^+\times (A^{(n-1)}\times{\cal P}(I))^+\ldots
\times (A^{(0)}\times{\cal P}(I))^+\times{\cal P}(I)$, where
$x
\in A^n$ is sent onto $(x,J_n(x),\ldots,J_0\circ J_1\circ\ldots J_n(x))$ and
$x\in Y^r$ is sent onto  $(\infty,\ldots,\infty,x,J_{r-1}(x),\ldots,
J_0\circ\ldots J_{r-1}(x)))$. This topological space will be denoted by
$\tilde X_n=A^{(n)}\sqcup_{J_n} Y_n\ldots\sqcup_{J_1} Y_1\sqcup_{J_0}{\cal J}$.

One deduces the following description of the spectrum $X$ of the Boolean algebra
${\cal B}=\cup_{n=0}^\infty{\cal B}_n$.
We define $X_{A,{\cal J}}$ as the set of terminal paths, where a terminal path is
either an infinite path $(i_0,i_1,\ldots)$  or a pair
$(i_0\ldots i_n,J)$ consisting of a finite path $i_0\ldots i_n$ and a set $J\in{\cal J}$
containing
$i_n$. The empty path $\emptyset$ is considered as finite path of length 0. The
corresponding terminal paths are $(\emptyset,J)$ with $J\in{\cal J}$. We note that
$X_{A,{\cal J}}$ is the projective limit of the $\tilde X_n$'s, where the projection
$\pi_{n,n+1}:\tilde X_{n+1}
\rightarrow \tilde X_n$
sends $(i_0,\ldots,i_n,i_{n+1})\in A^{(n+1)} $ and
$(i_0,\ldots,i_n;J)\in Y_{n+1}$  onto $(i_0,\ldots,i_n)\in A^n$ and  is the identity
map elsewhere and we endow it with the projective limit topology. We
define the symbol map $\sigma:X\rightarrow X_{A,{\cal J}}$ by looking at the orbit
$(x,Tx, T^2x,\ldots)$ of $x\in X$. We call exit time
$\tau=\tau(x)$ the smallest integer, if it exists, such that
$T^\tau x\notin U$; if it doesnot exist, we set $\tau(x)=\infty$. For $n<\tau(x)$, we
define
$i_n(x)=i_0(T^nx)$. If $x$ has infinite exit time, we define $\sigma(x)$ as the infinite path
$(i_0(x),i_1(x),\ldots)$ and if $x$ has a finite exit time
$\tau$, we define $\sigma(x)$ as the terminal path
$(i_0(x),\ldots),i_{\tau-1}(x);I(T^\tau x))$. We note that the symbol maps $\sigma_n$
define an isomorphism of the projective systems
$(X_n)$ and $(\tilde X_n)$ and that $\sigma$ is their limit. Therefore, it is a
homeomorphism. The symbol map $\sigma$ conjugates $T$ and the one-sided shift
$T_{A,{\cal J}}$ on
$X_{A,{\cal J}}$, defined on the set $\tilde U=\sigma(U)$ of terminal paths of length
$\ge 1$ by
$T_{A,{\cal J}}(i_0,\alpha)=\alpha$. 

We summarize the above discussion.

\begin{prop}\label{model}
Let $I$ be a countable set, $A$ a subset of $I\times I$ and
$\cal J$ a subset of
${\cal P}(I)$ containing the set ${\cal J}_A$ of cluster points of the net $(A_i=\{j:
(j,i)\in A\}),i\in I$. Assume that for each $i\in I$, either there exists $j\in I$ such that
$(i,j)\in A$ or there exists $J\in {\cal P}(I)$ such that $i\in J$. Then the SGDS
$(X_{A,{\cal J}},T_{A,{\cal J}})$ admits the Markov partition $\{\tilde U_i\},i\in I$,
where $\tilde U_i$ is the set of terminal paths starting with
$i$.
 
Moreover, any SGDS $(X,T)$ which admits a Markov
partition $\{ U_i\},i\in I$ is conjugate to the above model, where
$A=\{(j,i)\in I\times I: V_j\supset U_i\}$ and 
${\cal J}=\{J\in{\cal P}(I):
\cap_JV_j\cap_{I\setminus J}V^c_k\cap U^c\not=\emptyset\}$.
\end{prop}

We are chiefly interested in Markov shifts, i.e. SGDS $(X,T)$ admitting a Markov
partition and having a dense domain.

\begin{prop}\label{dense domain}
With the notation  of the above proposition, the domain
$U=\cup_I U_i$ of the SGDS
$(X_{A,{\cal J}},T_{A,{\cal J}})$ is dense iff ${\cal J}={\cal J}_A$.
\end{prop}

\begin{proof}  This is clear since, keeping the notation of the above discussion, $x\in
U^c$ is in the closure of $U$ iff $\sigma_0(x)\in J_A$.
\end{proof}

In this case, we shall
simply write
$X_A=X_{A,{\cal
J}}$ and $T_A=T_{A,{\cal J}}$. Note that in this case the assumption of the part $(i)$ of
the above proposition is simply that for each
$i\in I$,
$A^i=\{j\in I: (i,j)\in A\}$ is nonempty. 
When we view $I$ as the set of vertices and $A$ as the set of edges of an oriented
graph, this condition means that every vertex has at least one outgoing edge. When we
view $A$ as a matrix, this means that it has no zero rows.

The above analysis is also valid when $I$ is finite. In this case ${\cal J}_A$ is empty. If
one chooses ${\cal J }=\emptyset$, one gets the usual Markov shift $(X_A,T_A)$. If one
chooses ${\cal J }=\{I\}$, one gets the shift on the space of all sequences (finite or
infinite) which gives the Toeplitz extension.

Exel and Laca give a necessary and sufficient condition on the graph $A\subset I\times
I$ ensuring that the Markov shift $(X_A,T_A)$ is essentially principal
(\cite{el:infinitematrices}, Proposition~12.2). This condition also appears in
\cite{kpr:graphs} in the case of a row-finite directed graph, where it is called exit
condition (L).  A loop of length $n\ge 1$ is a finite path
$\alpha=(i_0,\ldots,i_n=i_0)$. An edge
$(i_k,j)$, where
$j\not= i_{k+1}$ (where $i_{k+n}=i_k$) is called an outgoing edge of the loop
$\alpha$. A point $x\in X$ is called periodic with respect to the SGDS $(X,T)$
if there exist $m<n$ such that $x$ belongs to the domain of $T^n$ and $T^mx=T^nx$. Then,
we say that $n-m$ is a period of $x$. The periodic
points of the Markov shift
$(X_A,T_A)$ are the infinite paths which,after some time, repeat indefinitely the same
loop $\alpha=(i_0,\ldots,i_p=i_0)$ and such an infinite path is
isolated if and only if the loop $\alpha$ has no outgoing edge. Thus,

\begin{prop}\label{condition[L]}
Let $(X_A,T_A)$ be the Markov shift constructed from the
graph 
$A\subset I\times I$ such that every vertex has at least one outgoing edge. The
following conditions are equivalent:
\begin{enumerate}
\item $(X_A,T_A)$  is essentially principal;
\item  condition (L): every loop has at least one outgoing edge;
\item $(X_A,T_A)$ has no isolated periodic point.
\end{enumerate}
\end{prop}

Thus, when the graph $A\subset I\times I$ satisfies these conditions,  
the groupoid of germs $Germ(X_A,T_A)$ coincide with the semidirect product 
$G(X_A,T_A)$.
 
Let us study next the relation between a set $\{S_i, i\in I\}$ of  partial
isometries on a Hilbert space $\cal H$ satisfying  the
Cuntz-Krieger relations associated to
$A\subset I\times I$ and the C$^*$-algebra $C^*(X_A,T_A)$. In order to avoid infinite
sums, Exel and Laca have introduced the following version of the Cuntz-Krieger relations,
where $P_i=S_iS_i^*, Q_i=S_i^*S_i$:
\begin{enumerate}
\renewcommand{\labelenumi}{(CK \arabic{enumi})}
\item the $Q_i$'s commute; 
\item the $P_i$'s are pairwise orthogonal;  
\item $P_jQ_i=A(i,j)P_j$ for all $(i,j)\in I\times I$, where $A$ is identified with
its characteristic function;
\item $\prod_{j\in E} Q_j\prod_{k\in F} (1-Q_k)=\sum_{i\in I} A(E,F,i)P_i$ for all
finite subsets
$E,F\subset I$ such that $A(E,F,i)=\prod_{j\in E} A(j,i)\prod_{k\in F} (1-A(k,i))$ is
nonzero except for a finite number of $i$'s.
\end{enumerate}

Let $(X,T)$ be a SGDS admitting the Markov partition $\{U_i, i\in I\}$. We define
$A\subset I\times I$ by $(i,j)\in A$ iff $U_j\subset V_i=T(U_i)$.
For each $i\in I$, we have the bisection of $G(X,T)$:
$$ S_i=\{(x,1,Tx), x\in U_i\}.$$
By definition, the ranges $U_i=r(S_i)$ and the domains $V_i=s(S_i)$ of the $S_i$'s
satisfy:
\begin{enumerate}
\setcounter{enumi}{1}
\renewcommand{\labelenumi}{(CK' \arabic{enumi})}
\item the $ U_i$ are pairwise disjoint;
\item for all $(i,j)\in I\times I$,  $U_j\subset V_i$ or $U_j\subset
V^c_i$ according to $A(i,j)=1$ or $A(i,j)=0$.
\end{enumerate}

The following lemma elucidates the meaning of the condition $(CK4)$.

\begin{lem} Let $(X,T)$ be a SGDS admitting a Markov partition $\{U_i, i\in I\}$.
Then, its domain $U=\cup_I U_i$ is dense iff
\begin{enumerate}
\setcounter{enumi}{3}
\renewcommand{\labelenumi}{(CK' \arabic{enumi})}
\item whenever the intersection $\cap_FV_j\cap_GV^c_k$, where $F,G$ are finite
subsets of $I$, contains finitely many $U_i$'s, it is the union of these $U_i$'s.
\end{enumerate}
\end{lem}

\begin{proof} Suppose that $U$ is dense. If the intersection $\cap_FV_j\cap_GV^c_k$,
where
$F,G$ are finite subsets of $I$, contains only a finite union $\cap_EU_i$, then the
complement of
$\cap_EU_i$ in $\cap_FV_j\cap_GV^c_k$, which is open, must be empty. On the other
hand, suppose that $(iv')$ holds. We shall show that
${\cal J}={\cal J}_A$, which by \propref{dense domain} is equivalent to the density of
$U$. If
$J\in {\cal P} (I)$ is not a cluster point of the map $i\mapsto J(i)$ introduced in the
discussion following \defnref{markov}, there exist finite subsets $F,G$ of $I$ such that
$F\subset J, G\subset J^c$ and the conditions
$F\subset J(i), G\subset J(i)^c$ are satisfied for only a finite number of $i$'s. This
means that $\cap_FV_j\cap_GV^c_k$ contains only finitely many $U_i$'s; by $(iv')$, it
is contained in $U$. By definition, $J$ does not belong to ${\cal J}$.
\end{proof}

\begin{prop} 
Let $I$ be a countable set and let $A$ be a subset of $I\times
I$.
\begin{enumerate}
\item Let $\cal J$ be a subset of
${\cal P}(I)$ containing the set ${\cal J}_A$.  Then the family of
bisections $(S_i), i\in I$ of
$G(X_{A,{\cal J}},T_{A,{\cal J}})$, viewed as elements of 
$C^*(X_{A,{\cal J}},T_{A,{\cal J}})$ is a family of partial
isometries which  satisfies the Cuntz-Krieger relations
$(CK1-3)$ relative to the matrix $A$. In particular, every representation
of
$C^*(X_{A,{\cal J}},T_{A,{\cal J}})$ on a Hilbert space $\cal H$ provides a representation
of these relations.
\item Conversely, let $I$ be a countable set, let $A$ be a $(0,1)$-matrix with
no zero rows and let $(\underline S_i), i\in I$ be a family of nonzero partial isometries
on a Hilbert  space $\cal H$ satisfying the Cuntz-Krieger relations $(CK1-3)$ relative to
the matrix
$A$. Then, there exists a unique representation of $C^*(X_{A,{\cal J}},T_{A,{\cal J}})$,
where
$${\cal J}=\{J\in{\cal P}(I):E,F \hbox{finite}\  E\subset J, F\subset J^c\Rightarrow
\prod_E Q_j\prod_F(1-Q_j)\not=0\}.$$ on 
$\cal H$ sending the bisection $S_i$ into  $\underline S_i$ for each $i$.
\end{enumerate}
\end{prop}

\begin{proof} (Cf. Proposition~III.2.7 of \cite{ren:approach} and Theorem~4.2 of
\cite{kprr:graphs}) For
$(i)$, the condition $(CK1)$ is satisfied by construction. The conditions $(CK2-3)$ are a
restatement of above properties $(CK'2-3)$.

For $(ii)$, we give the sketch of the proof and refer to  \cite{ren:approach},
\cite{el:infinitematrices} and
\cite{kprr:graphs} for details. Let $(\underline S_i), i\in I$ be such a family of
partial isometries. A direct computation (see Proposition 3.2 of
\cite{el:infinitematrices}) shows that the
$\underline S_i$'s (together with 1) generate an inverse semigroup 
$\underline{\cal S}$
of partial isometries. Let
$X$ be the spectrum of the (commutative) C$^*$-algebra
$\cal B$ generated by the idempotents of $\cal S$ and let $M:C(X)\rightarrow
{\cal B}$ be the Gelfand isomorphism. For each
$i$, let
$U_i\subset X$ [resp. $V_i\subset X$] be the support of $\underline P_i$ [resp.
$\underline Q_i$]. The isomorphism
$B\mapsto\underline S_iB\underline S_i^*$ from $\underline Q_i{\cal B}$ onto
$\underline P_i{\cal B}$ induces a homeomorphism
$T_i:U_i\rightarrow V_i$ and we define $T$ on the (disjoint) union $U=\cup_I U_i$ by
$Tx=T_ix$ if $x\in U_i$. Then $\{U_i, i\in I\}$ is a Markov partition for $(X,T)$ with
matrix $A$ . According to \propref{model},  $(X,T)$ is isomorphic to
$(X_{A,{\cal J}},T_{A,{\cal J}})$. Let
$\{S_i, i\in I\}$ be the associated family of bisections of $G(X,T)$. It satisfies the same
relations $(CK1-3)$ as the family
$\{\underline S_i, i\in I\}$. Let $\cal S$ be the inverse semigroup of 
bisections of $G(X,T)$ generated by the $S_i$'s (and $X$). Given two finite paths
$\alpha=(i_1,\ldots,i_m)$ and
$\beta=(j_1,\ldots,j_n)$ and  an idempotent $h$ of $\cal S$, we define
$$S(\alpha,h,\beta)=S_{i_m}^{-1}\ldots S_{i_1}^{-1}hS_{j_1}\ldots S_{j_n}\quad{\rm
and}$$
$$\underline S(\alpha,h,\beta)=\underline S_{i_m}^*\ldots
\underline S_{i_1}^*M(h)\underline S_{j_1}\ldots \underline S_{j_n}.$$
One checks that every element of $\cal S$ can be written under the form
$S(\alpha,h,\beta)$ and that this form is unique when we require $h$ and the lengths
of $\alpha,\beta$ to be minimal. The same is true for $\underline{\cal S}$. Therefore
the map
$S(\alpha,h,\beta)\mapsto
\underline S(\alpha,h,\beta)$ is a representation
$L$ of the inverse semigroup $\cal S$ onto
$\underline{\cal S}$. This representation extends to a representation of $C^*(G(T))$.
Indeed, every $f\in C_c(G(T))$ can be written as a finite sum $f=
\sum (h_\alpha\circ r) S_\alpha$, where $h_\alpha\in C(X)$ and $S_\alpha$ is a
compact open bisection. One can check that $L(f)=\sum M(h_\alpha)L(S_\alpha)$ is well
defined and that $L$ is a representation of the $*$-algebra $C_c(G(T))$ continuous for
the inductive limit topology and therefore extends to $C^*(G(T))$.
\end{proof}

\begin{cor} Let $I$ be a countable set and let $A$ be a subset of $I\times I$.
\begin{enumerate}
\item The family of
bisections $(S_i), i\in I$ of
$G(X_A,T_A)$, viewed as elements of 
$C^*(X_A,T_A)$ is a a family of partial
isometries which  satisfies the Cuntz-Krieger relations
$(CK1-4)$ relative to the matrix $A$. In particular, every representation
of
$C^*(X_A,T_A)$ on a Hilbert space $\cal H$ provides a representation
of these relations.
\item Conversely, let $I$ be a countable set, let $A$ be a $(0,1)$-matrix with
no zero rows and let $(\underline S_i), i\in I$ be a family of nonzero partial isometries
on a Hilbert  space $\cal H$ satisfying the Cuntz-Krieger relations $(CK1-4)$ relative to
the matrix
$A$. Then, there exists a unique representation of $C^*(X_A,T_A)$,
on  $\cal H$ sending the bisection $S_i$ into  $\underline S_i$ for each $i$.
\end{enumerate}
\end{cor}

\begin{proof}  We have seen indeed the equivalence of the conditions $(CK4)$ and
${\cal J}= {\cal J}_A$.
\end{proof}

Because of this corollary, the C$^*$-algebra $C^*(X_A,T_A)$, where $A\subset I\times
I$ has no zero rows, can be called the Cuntz-Krieger algebra of the  matrix
$A$. If moreover $A$ the equivalent conditions of Proposition~\ref{condition[L]}, we
have the following uniqueness result.

\begin{cor} (Cf. \cite{el:infinitematrices}, 13.2) Let $A\subset I\times
I$ without zero rows and such that each loop has an outgoing edge. Let
$\{ S_i), i\in I\}$ and $\{ T_i), i\in I\}$  be two families of nonzero partial isometries
on  Hilbert  spaces $\cal H$ and $\cal K$ satisfying the Cuntz-Krieger relations
$(CK1-4)$ relative to the matrix
$A$. Then there is an isomorphism from the C$^*$-algebra generated by $\{ S_i), i\in
I\}$ onto the C$^*$-algebra generated by 
$\{ T_i), i\in I\}$ carrying $S_i$ into $T_i$ for all $i\in I$.
\end{cor}

\begin{proof} It suffices to show that the representation $L$ of $C^*(X_A,T_A)$
defined by
$\{ S_i), i\in I\}$ is faithful. Since its restriction $M$ to $C(X_A)$ is faithful and
the groupoid $G(X_A,T_A)$ is nuclear and essentially principal, we may apply
Theorem~4.3 of \cite{ren:ideal} to conclude.
\end{proof}

 By construction,
$C^*(X_A,T_A)$ is unital and generated by the partial isometries $S_i$'s and $1$. Exel
and Laca distinguish this unital C$^*$-algebra, which they call $\tilde O_A$ and the
sub-C$^*$-algebra $O_A$ generated by the
$S_i$'s alone. There are two cases: either $\tilde O_A= O_A$ or $\tilde O_A$ is the
one-point compactification of $O_A$. Recall that in the description of the spectrum
$X=X_A$, we have introduced the map $\sigma_0:X\rightarrow X_0=I\sqcup{\cal J}$,
which is a homeomorphism from $U^c$ onto ${\cal J}$. Its inverse is the map $J\mapsto
(\emptyset;J)$ where the first component means the empty word.  The complement of
$U\cup V$ contains at most one point, namely the point
$(\emptyset;\emptyset)$. The spectrum of the C$^*$-algebra generated by the $P_i$'s and
the
$Q_i$'s alone is either $X_A$ or $X_A\setminus\{(\emptyset;\emptyset)\}$ according to
$\emptyset\notin {\cal J}_A$ or $\emptyset\in {\cal J}_A$. Thus we
have:

\begin{prop}(Cf. \cite{el:infinitematrices}, 8.5) Let $A$ be a $(0,1)$-matrix with
no zero rows.
\begin{enumerate}
\item if $\emptyset\notin {\cal J}_A$, then $O_A=C^*(X_A,T_A)$; 
\item if $\emptyset\in {\cal J}_A$, then
$O_A=C^*(X_A\setminus\{(\emptyset;\emptyset)\},T_A)$.
\end{enumerate}
\end{prop}

We may apply to the SGDS $(X_A,T_A)$ the general criteria for simplicity and pure
infiniteness given in the first section.

\begin{prop}(Cf. \cite{el:infinitematrices},14.1) Assume that the matrix $A\subset
I\times I$ satisfies
$(L)$ and is irreducible, in the sense that for each pair $(i,j)\in I\times I$, there is a
finite path starting at $i$ and ending at $j$. Then $O_A$ is simple.
\end{prop}

\begin{proof} The SGDS $(X_A,T_A)$ and
$(X_A\setminus\{(\emptyset,\emptyset)\},T_A)$ are essentially free and we
may apply \propref{simplicity}.  Every nonempty open set
$W$ contains a cylinder set, i.e. a set $Z(\alpha)$ of all terminal paths starting with the
finite path
$\alpha=(j_0,\ldots,j_n)$. Let
$x\in X_A$ be an infinite path $x=(i_0,i_1,\ldots)$. There exist a finite  path 
$\beta=(j_n, \ldots,j_{n+k}=i_0)$. Then
$y=(j_0,\ldots,j_n,\ldots,j_{n+k}=i_0,i_1,\ldots)$ belongs to $Z(\alpha)$, hence to $W$,
and
$T^{n+k}y=x$. Let
$x$ be a finite terminal path distinct from $(\emptyset;\emptyset)$. There exists an
integer $m$ such that $T^mx=(\emptyset;J)$ where $J$ is nonempty. We choose
$i\in J$ and a finite  path 
$\beta=(j_n, \ldots,j_{n+k}=i)$. Then
$y=(j_0,\ldots,j_n,\ldots,j_{n+k}=i;J)$ belongs to  $W$
and
$T^{n+k}y=T^mx$.
\end{proof}

\begin{prop}(Cf. \cite{kpr:graphs}, 3.9  and \cite{el:infinitematrices},16.2) Assume that
the matrix
$A\subset I\times I$ satisfies
$(L)$ and for every vertex $i\in I$, there exists a loop $(j_0,\ldots,j_n=j_0)$ and a finite
path starting at
$i$ and ending at some $j_k$. Then
$C^*(X_A,T_A)$ and $O_A$ are purely infinite.
\end{prop}

\begin{proof} We apply \propref{pure infiniteness} to the SGDS $(X_A,T_A)$. First observe
that if the loop
$(j_0,\ldots,j_n=j_0)$ has an outgoing edge, then 
$Z(j_0,\ldots,j_n)$ is a
proper subset of $T^nZ(j_0,\ldots,j_n)$. We have seen that every nonempty open set $W$
contains a cylinder set
$Z(i_0,\ldots,i_m)$. By assumption, there exists a loop $(j_0,\ldots,j_n=j_0)$ and a
finite path
$i_m,\ldots,i_{m+k}=j_0$. Let $V=Z(i_0,\ldots,i_{m+k}=j_0,\ldots,j_n)$. Then $V\subset
W$ and $T^{m+k}V$ is a proper subset of $T^{m+k+n}V$.
\end{proof}


\begin{thebibliography}{10}

\bibitem{del:purelyinfinite}  C.~ Anantharaman-Delaroche: {\it Purely infinite
$C^*$-algebras arising
from dynamical systems}, Bull. Soc. Math. France, {\bf 125} (1997), 199-225.

\bibitem{dr:amenable}  C.~ Anantharaman-Delaroche and J.~Renault: {\it Amenable
groupoids}, preprint 1998.

\bibitem{ar:examples}  V.~ Arzumanian and J.~Renault: {\it Examples of pseudogroups
and their $C^*$-algebras}, in Operator Algebras and Quantum Field Theory, 
S.~Doplicher, R.~Longo, J.~E.~Roberts and L.~Zsido, editors, International Press
1997, 93-104.

\bibitem{cfw:amenable}  A.~Connes, J.~Feldman and B.~Weiss: {\it An amenable
equivalence relation is generated by a single transformation}, J.~Ergodic Theory and
Dynamical Systems, {\bf 1} (1981), 431-450.


\bibitem{cun:simple}  J.~Cuntz: {\it Simple $C^*$-algebras generated by
 isometries}, Comm. Math. Phys., {\bf  57} (1977), 173-185.

\bibitem{cun:markovII}  J.~Cuntz: {\it A class of $C^*$-algebras and
topological Markov chains II: reducible chains and the Ext-functor for $C^*$-algebras},
Inventiones Math., {\bf  63} (1981), 25-40.

\bibitem{cun:homotopy}  J.~Cuntz: {\it On the homotopy groups of the space of
endomorphisms of a $C^*$-algebra (with applications to topological Markov chains)} in
{\it Operator algebras and group
representations}, Proc. Int. Conf., Vol.~{\bf I} Pitman Boston
(1984), 124-137.

\bibitem{ckr:markovI}  J.~Cuntz and W.~Krieger: {\it A class of $C^*$-algebras and
topological Markov chains},
Inventiones Math., {\bf  56} (1980), 251-268.

\bibitem{dea:groupoids} V.~Deaconu: {\it Groupoids associated with
endomorphisms}, Trans. Amer. Math. Soc., {\bf 347} (1995), 1779-1786.

\bibitem{el:infinitematrices}  R.~Exel and M.~Laca: {\it Cuntz-Krieger algebras
for infinite matrices}, preprint 1998.

\bibitem{fm:relations}  J.~Feldman and C.~Moore: {\it Ergodic equivalence relations,
cohomologies, von Neumann algebras, I and II},
Trans. Amer. Math. Soc., {\bf 234}, no. 2 (1977), 289-359.

\bibitem{hae:pseudogroups}  A.~Haefliger: {\it Pseudogroups of local isometries},
Pitman Research Notes  {\bf 131} (1985), 174-197.

\bibitem{kri:dimensions}  W.~Krieger: {\it On dimension functions and topological
Markov chains},  Inventiones Math., {\bf  56} (1980), 239-250.

\bibitem{kum:localizations}  A.~Kumjian: {\it On localizations and simple
C*-algebras}, Pacific J. Math. {\bf 112} (1984), 141-192.

\bibitem{kprr:graphs}  A.~Kumjian, D.~Pask, I.~Raeburn and J.~Renault: {\it Graphs,
groupoids and Cuntz-Krieger algebras}, J. Funct. Anal.
{\bf 144} (1997), 505-541.

\bibitem{kpr:graphs}  A.~Kumjian, D.~Pask and I.~Raeburn: {\it Cuntz-Krieger
algebras and directed graphs}, Pacific J.~Math., to appear.

\bibitem{lm:intro}  D.~Lind and B.~Marcus: {\it An introduction to symbolic
dynamics and coding}, Cambridge University Press, Cambridge, 1995.


\bibitem{mat:subshifts}  K.~Matsumoto: {\it  On $C^*$-algebras associated with
subshifts}, International J.~Math. {\bf 8} No.3 (1997), 357-374.


\bibitem{moo:relations}  C.~C.~Moore: {\it Equivalence relations and von Neumann
algebras} in {\it Operator algebras and applications }, Proc.~Sympos. Pure Math.,
Vol.~{\bf 38}, Part ~{\bf II}, American Mathe\-matical Society R.I.
(1982), 339-350.

\bibitem{mrw:equivalence}  P.~Muhly, J.~Renault and D.~Williams: {\it Equivalence
and isomorphism for groupoid $C^*$-algebras}, J.~Operator Theory {\bf 17}(1987), 3-22

\bibitem{ren:approach}  J.~Renault: {\it A groupoid approach to
$C^*$-algebras}, Lecture Notes in Mathematics, Vol.~{\bf 793}
Springer-Verlag Berlin, Heidelberg, New York (1980).

\bibitem{ren:ideal}  J.~Renault: {\it  The ideal structure of groupoid crossed
product
$C^*$-algebras}, J.~Operator Theory {\bf 25}
(1991), 3-36.

\bibitem{rie:rotations}  M.~Rieffel: {\it $C^*$-algebras associated with irrational
rotations}, Pacific J.~Math. {\bf 98}, No 2 (1981), 415-429.




\end{thebibliography}
\end{document}